\documentclass[10pt]{article}

 \usepackage[utf8]{inputenc}
\usepackage[english]{babel}
\usepackage{amsmath, amsfonts, amstext, amssymb} 

 \usepackage{makeidx,epsfig,lscape}
 \usepackage{color,colortbl}
 \usepackage{fancyhdr}
 \usepackage{xcolor,pict2e}

 \usepackage{tikz}
\usepackage[hypcap]{caption}
\usetikzlibrary{shapes}
\usepackage{lipsum}
\usepackage[all,dvips]{xy}
\usepackage{graphicx}
 
 \renewcommand{\headrulewidth}{0pt}
 \renewcommand{\footrulewidth}{0.5pt}

 \definecolor{myaqua}{rgb}{0.0,0.5,0.55}
 \definecolor{lightaqua}{rgb}{0.75,0.95,0.95}

 \usepackage[colorlinks = true,
            linkcolor = myaqua,
            urlcolor  = blue,
            citecolor = red,
            pdfpagemode=UseOutlines,
            bookmarksnumbered=true,
            pdfpagelabels,
            breaklinks]{hyperref}

\usepackage{caption}
\usepackage{floatrow}

 \newcommand{\N}{\mathbb{N}}

\newcommand{\R}{\mathbb{R}}

\newcommand{\A}{\mathcal{A}}

\newtheorem{theo}{Theorem}
\newtheorem{coro}{Corollary}
\newtheorem{lem}{Lemma}
\newtheorem{defi}{Definition}[section]
\newtheorem{prop}{Proposition}[section]

\newtheorem{rem}{Remark}[section]

\newenvironment{preu}{\textbf{Proof }\rm}{$\blacksquare $}
\newenvironment{preu 2}{\textbf{Proof for $ d_{2i}= d_1 $.}\rm}{$\blacksquare $}
\newenvironment{preu 3}{\textbf{Proof of theorem 1}\rm}{$\blacksquare $}
\newenvironment{preu 4}{\textbf{Proof of lemma 1}\rm}{$\blacksquare $}
\newtheorem{exple}{Example}

\def\lin#1#2{\textcolor[rgb]{0.6,0.6,0.6}{\vspace*{#1mm} \hrule
   height 3 pt \vspace*{#2mm}}}
\def\bt{\begin{tabular}}
\def\et{\end{tabular}}
\def\and{\mbox{ and }}

\def\1{{\bf 1}}

 \def\boxx#1#2#3#4#5{
 {\linethickness{#4pt}\put(#1,#5){\color{myaqua}{\line(1,0){#3}}}}
 \multiput(#1,#2)(0,#4){2}{\line(1,0){#3}}
 \multiput(#1,#2)(#3,0){2}{\line(0,1){#4}}
  }

\begin{document}


 $\mbox{ }$

 \vskip 12mm

{ 

{\noindent{\Large\bf\color{myaqua}
   On the Lyndon dynamical system  }} 
%
\\[6mm]
{\bf Florent NGUEMA NDONG}}
\\[2mm]
{ 
  Université des Sciences et Techniques de Masuku 
 \\
Email: \href{mailto:florentnn@yahoo.fr}{\color{blue}{\underline{\smash{florentnn@yahoo.fr}}}}\\[1mm]
\lin{5}{7}

 {  
 {\noindent{\large\bf\color{myaqua} Abstract}{\bf \\[3mm]
 \textup{
Given a totally finite ordered alphabet $ \A $, endowing the set of words over $ \A $ with the alternating lexicographic order, we define a new class of Lyndon words. We study the fundamental properties of the associated symbolic dynamical systems called Lyndon system. We derive some fundamental properties of the beta-shift with negative base by relating it with the Lyndon system. We find, independently of W. Steiner's method, the conditions for which a word is the $(-\beta)$-expansion of $-\frac{\beta}{\beta + 1} $ for some $ \beta> 1$.
 }}}
 \\[4mm]
 {\noindent{\large\bf\color{myaqua} Keywords}{\bf \\[3mm]
 Lyndon Word; Lyndon System; Dynamical System; Symbolic dynamics; expansion in negative base
}}}
\lin{3}{1}

\renewcommand{\headrulewidth}{0.5pt}
\renewcommand{\footrulewidth}{0pt}

 \pagestyle{fancy}
 \fancyfoot{}
 \fancyhead{} 
 \fancyhf{}
 \fancyhead[RO]{\leavevmode \put(-90,0){\color{myaqua}F. NGUEMA NDONG} \boxx{15}{-10}{10}{50}{15} }
 \fancyfoot[C]{\leavevmode
 \put(-2.5,-3){\color{myaqua}\thepage}}

 \renewcommand{\headrule}{\hbox to\headwidth{\color{myaqua}\leaders\hrule height \headrulewidth\hfill}}
 
\section{Introduction}

In the areas of combinatorics and computer sciences, a Lyndon word is a word which is lexicographically less than all its permutations. Roger Lyndon introduced them in 1954 on the standard name ``lexicographic sequences''. He used them to construct a basis for the homogeneous part of a given degree in Lie algebra (see \cite{MR0064049,MR0102539}). 

By definition, we call Lyndon word $x_1x_2 x_3 \cdots $ over a totally ordered alphabet all word which is lexicographically less than all of its suffixes. 

\begin{equation*}
x_1 x_2 \cdots x_{n-1}x_n \cdots \leq_{lex} x_k x_{k+1} \cdots x_nx_{n +1} \cdots.
\end{equation*}

Consider a finite ordered alphabet $ \A = \{ 0, 1, \cdots, d \} $. We note $ \A^{\N} $ the set of finite and infinite words over $ \A $. Endowing $\A^{\N}$ with the alternate order on words (see \cite{MR2223030,MR2534912}), we can define Lyndon word using this new order. We call them alternate Lyndon words. The definition given in \cite{MR2223030} generalizes Lyndon words. To each Lyndon word we attach a symbolic dynamical system. If the entropy of such a system is positive (shall we say $ \log \beta>0 $), we show that the associated alternate Lyndon word is a $(-\beta)$-representation of $-\dfrac{\beta}{\beta + 1}$ (proposition \ref{proposition 3}). When $ \beta $ tends to 1, the alternate Lyndon word tends to 
\begin{equation*}
\phi^{\infty}(1) = 1001110010010011100111001110010010011100100\cdots,
\end{equation*}
where $ \phi $ is a morphism on $\{0, 1\}$ such that $ \phi(0) = 1 $ and $ \phi(1) = 100 $ (theorem \ref{theorem 1}). We establish a link with expansions in negative bases. Indeed, For a fixed real $ \beta > 1$, we can see the alternate order as a tool of controllability of the representations of numbers in base $ -\beta $. In the negative expansion pioneering paper \cite{MR2534912}, S. Ito and T. Sadahiro proved that for some $ \beta > 1$, the $(-\beta)$-expansion of $ -\frac{\beta}{\beta + 1} $,  $ (d_i)_{i \geq 1}$ is such that all sub-word of a $(-\beta)$-expansion $ (x_i)_{i \geq 1}$ is greater than $ (d_i)_{i \geq 1}$ in the sense of the alternate order. In particular, $(d_i)_{i \geq 1}$ is less than all of its sub-words (in the sense of the alternate order). Thanks to this link, we give the necessary and sufficient conditions for which a word over a finite totally ordered alphabet can be the $(-\beta)$-expansion of $ - \frac{\beta}{\beta + 1}$ for some $\beta> 1$ (theorem \ref{theorem 3}).  

\subsection{Definitions and generality}

\begin{defi}

Let $ \A $ be a totally ordered alphabet endowed with an order "$<$". 
We call \textit{ Lyndon word} over $ \A $ with respect to the order "$<$", all word $ x_1 x_2 x_3 \cdots $ such that
\begin{equation}
x_1 x_2 x_3 \cdots \leq x_n x_{n+1} x_{n+2} \cdots,  \text{ $ \forall n $ }. 
\label{1}
\end{equation}
with $ x_i $ in $ \A$ for all $i$.
\begin{itemize}
\item The word $ x_1 x_2 x_3 \cdots $ is said \textit{strong Lyndon word } if in \eqref{1} all inequalities are strict.
\item The Lyndon word $ x_1 x_2 \cdots $ is weak if in \eqref{1} equality holds for some $ n $. 
\end{itemize}

\begin{equation*}
 x_1x_2 x_3 x_4 \cdots = x_n x_{n+1} x_{n+2} \cdots.
\end{equation*} 
\end{defi}
Then, the  weak Lyndon words over the alphabet $ \A $  are periodic.

\subsection{Ito and Sadahiro order}

Let $ \A = \{0, 1, \cdots, d \} $ be an alphabet. Consider two words $ x_1 x_2 \cdots x_n $ and $ y_1 y_2 \cdots y_n $ on $ \A $. We will say $ x_1 x_2 \cdots x_n $ is less than $ y_1 y_2 \cdots y_n $ in the sense of alternate order (and we note $ x_1 x_2 \cdots x_n \prec y_1 y_2 \cdots y_n $) if there exists an integer $ k \leq n $ such that for all $ i < k $, $ x_i = y_i $ and $ (-1)^k(x_k-y_k)< 0$. We will note $ x_1 \cdots x_n \preceq y_1 \cdots y_n $ if $ x_1 \cdots x_k = y_1 \cdots y_n $ or $ x_1 \cdots x_n \prec y_1 \cdots y_n$. 
The relation "$ \preceq$" is called \textit{Ito and Sadahiro order} or \textit{alternate order}.

We can extend the alternating order to infinite words and it is possible to compare two words with different lengths by completing the least at right by a sequence of zeros. Indeed, consider $ x_1 \cdots x_k $ and $ y_1 y_2 \cdots y_n $ over $ \A $, $ k \neq n $. 
\begin{equation}
x_1\cdots x_k\preceq y_1\cdots y_n\Leftrightarrow \begin{cases}
                                                 x_1 \cdots x_k (0)^{n-k} \preceq y_1 \cdots y_n &\text{ if $ k< n $ } \\
                                                 x_1 \cdots x_k  \preceq y_1 \cdots y_n(0)^{k-n} &\text{ if $ n < k $}.
                                                  \end{cases}
\end{equation}

\begin{defi}
Let $ (d_i)_{i \geq 1} $ be a Lyndon word for an order $ \leq $ over $ \A $ (weak or strong). We call \textit{dynamical Lyndon system} associated to $(d_i)_{i \geq 1}$ with respect to the order $\leq $, the set of infinite words $ x_1 x_2 \cdots $ on $ \A $  such that : 
 
 \begin{equation*}
   d_1 d_2 d_3 \cdots \leq x_k x_{k+1} \cdots, \text{ $ \forall k $ }.
   \end{equation*}

\end{defi}
Such a one-side system is non empty, invariant by the shift $\sigma: (x_i)_{i \geq k} \mapsto (x_{i + 1})_{i \geq k}$.

\begin{exple}
Let $ \A = \{0, 1, \cdots, d\} $ be an alphabet endowed with the usual order and $ \A^{\N} $ equipped with the alternate order. Let $(x_i)_{i \geq 1}$ be a sequence on $ \A^{\N} $ satisfying for all integer $ k $ :
 
 \begin{equation}
  x_1 x_2 x_3 \cdots \preceq x_k x_{k+1} \cdots.
  \label{(2)}
 \end{equation}
 
 Such a word is called, in the following sentences, Lyndon word with respect to alternate order or just \textbf{alternate Lyndon word}. An \textbf{alternate Lyndon system} associated to the alternate Lyndon word $ (d_i)_{i \geq 1} $ is the set of words $ (x_i)_{i \geq 1}$ such that 

 \begin{equation*}
 (d_i)_{i \geq 1} \preceq x_n x_{n+1} \cdots, \text{ $ \forall n $ }. 
 \end{equation*}
 \end{exple}

\begin{rem}\label{rem 1}
If $ (d_i)_{i \geq 1}$ is an alternate Lyndon word, then  for all integer $ k $, $ d_k \leq d_1 $. The language of the associated Lyndon system can be taken over the finite alphabet $ \{ 0, 1, \cdots, d_1 \}$.
\end{rem}

\begin{rem}\label{rem 2}

Let $ (d_i)_{i \geq 1} $ be an alternate Lyndon word over a finite alphabet $ \{0, 1, \cdots, d_1\} $ endowed with the usual order. The following propositions are equivalent :
 
 \begin{itemize}
  \item For all sequence  $ (x_i)_{i \geq 1}$ on $ \A$, and for all $ k \geq 1 $,
  \begin{equation*}
 d_1 d_2 d_3 \cdots \preceq x_k x_{k+1} \cdots .
\end{equation*}
 
  \item for all $ k \geq 1 $
  
\begin{equation}
 d_1 d_2 d_3 \cdots \preceq x_k x_{k+1} \cdots \preceq 0d_1d_2\cdots .
 \label{(3)}
\end{equation}
 \end{itemize}

\end{rem}

\begin{exple}
In this example, we focus on the case of the $\beta$-shift. The alphabet $ \A = \{ 0, 1, \cdots, \lfloor \beta \rfloor \} $ is endowed with the usual order and $ (a_i)_{i \geq 1} $, the $ \beta$-expansion of 1. It is well-known that for all integer $ k $ :
  \begin{equation*}
  a_ka_{k+1} \cdots \leq_{lex} a_1a_2 a_3 \cdots 
 \end{equation*}
It is a Lyndon word if we consider the order $ \leq_L $ :
\begin{equation*}
 (x_n)_{n\geq1} \leq_{L} (y_n)_{n \geq 1} \Leftrightarrow (y_n)_{n \geq 1} \leq_{lex} (x_n)_{n \geq 1} 
\end{equation*}

To each number $ \beta>1 $ corresponds a unique infinite Lyndon word and then a unique Lyndon system
\begin{equation*}
 X_{\beta} = \{ x_1 x_2 \cdots; \forall k \text{ $ x_1x_2 \cdots \leq_{lex} a_1a_2 \cdots $ } \}
\end{equation*}
except if $ \beta $ is a simple Parry number, that is, the sequence $ (a_i)_{i \geq 1}$ ends by zeros.
\begin{equation*}
 1 = \frac{a_1}{\beta} + \frac{a_2}{\beta^2} + \cdots + \frac{a_k}{\beta^k}.
\end{equation*}
 Since
\begin{equation*}
 1 = \frac{a_1}{\beta} + \cdots + \frac{a_{k-1}}{\beta^{k-1}} + \frac{a_k - 1}{\beta^k} + \frac{a_1}{\beta^{k+1}} + \cdots + \frac{a_{k-1}}{\beta^{2k-1}} 
 + \frac{a_k-1}{\beta^{2k}} + \cdots,
\end{equation*}
$ \beta$ is attached to  $ a_1a_2 \cdots a_k \overline{0}$ and
$ \overline{a_1 a_2 \cdots a_{k-1} (a_k-1)} $. The corresponding systems are slotted, that generated by the weak Lyndon word $ \overline{a_1 a_2 \cdots a_{k-1} (a_k-1)} $ being the smallest. This one defines the $ \beta$-shift. 
\end{exple}

There are a monotonic bijection between the strong Lyndon words with respect to the lexicographic order and the real numbers strictly greater than 1. On the other hand, each simple  Parry number corresponds (only) to two Lyndon words with respect to lexicographic order : a weak Lyndon and a strong Lyndon word. The exponential of the entropy of the system is the number associated to the Lyndon word.

Let $ \A $ be a finite or countable alphabet. The set of Lyndon words over $ \A $ with respect to the alternate order is totally ordered.

Consider two alternate Lyndon words $ (a_i)_{i \geq 1} $ and $ (b_i)_{i \geq 1}$. We suppose $ (a_i)_{i \geq 1}\neq(b_i)_{i \geq 1}$. Then, there exists an integer $ k $ such that
\begin{equation*}
 a_i = b_i  \text{ for all $ i < k $ and } (-1)^k(a_k - b_k) \neq 0 .
\end{equation*}
We have either $(-1)^k(a_k-b_k) < 0$ and $ (a_i)_{i \geq 1}\prec(b_i)_{i \geq 1}$ or $(-1)^k(a_k-b_k)>0$ and $ (b_i)_{i \geq 1}\prec(a_i)_{i \geq 1}$. 

The Reutenauer approach of the alternate order is a little bit different. Indeed, consider 
two words $ (x_i)_{i \geq 1} $ and $ (y_i)_{i \geq 1} $ over an alphabet $ \A $. Then, $ (x_i)_{i \geq 1} \prec (y_i)_{i \geq 1} $ if only if: 

\begin{equation*}
x_1 + \dfrac{1}{x_2 + \dfrac{1}{x_3 + \dfrac{1}{\ddots}}} < y_1 + \dfrac{1}{y_2 + \dfrac{1}{y_3 + \dfrac{1}{\ddots}}}.
\end{equation*}
 The alternating lexicographical order considered by Steiner \cite{MR3028656}  is the inverse of the alternate order considered in the present paper.
\begin{rem}\label{rem 3}
Let $ U $, $V_1 $, $V_2$ three words over a finite alphabet such that $ \vert U \vert <+\infty $. Then, 
\begin{equation}
U V_1 \prec U V_2 \Longleftrightarrow \begin{cases} 
                                      V_1 \prec V_2 &\text{ if $ \vert U \vert $ is even }\\
                                      V_2 \prec V_1 &\text{ if $ \vert U \vert $ is odd  }
                                      \end{cases} \label{111}
\end{equation}
\end{rem}

\section{Alternate Lyndon system}
\subsection{Generality}
Let $ (d_i)_{i \geq 1} $ be a Lyndon word with respect to the alternate order. We note $M=M((d_i)_{i \geq1 })$ the associated Lyndon system and $ L_M $ its language. That is, $ L_M $ denotes the set of finite sequences of words which belong to $ M=M((d_i)_{i \geq1 }) $. All word $ x_1 x_2 \cdots x_n $ of $ L_M $ satisfies :

\begin{equation}
 d_1 d_2 \cdots d_{n-j+1} \preceq x_j x_{j+1} \cdots x_n \preceq 0d_1d_2\cdots d_{n-j}, \text{ with $ 1 \leq j \leq  n $ }. 
 \label{(4)}
\end{equation}

The relation above implies that $ (x_i)_{ i\geq 1} $ is bounded (see in the remark \ref{rem 1}). The words of $ M $ are taken in the finite alphabet $\A$. Furthermore, the relation \eqref{(4)} implies that all sub-word of a word of $ M $ belongs to $ M $. That is, $M$ is invariant by the shift $ \sigma : (x_i)_{i \geq 1} \mapsto (x_{i+1})_{i\geq1} $. Since $ M $ is a closed subset of $ \A^{\N} $ endowed with the metric 
\begin{equation*}
d(u,w)= \sum \limits_{ n \in \N} \frac{d_n(u,w)}{2^{\vert n \vert}}, \text{ with } 
\end{equation*}

\begin{equation*}
 d_n(u, w) = \begin{cases}
                0 \text{ if $ u_n = w_n $ }\\
                1 \text{ if $ u_n \neq w_n $ },
               \end{cases}
\end{equation*}
it follows that $ (M, \sigma)$ (we often denote $ M $ if there is no ambiguity) is a symbolic dynamical system.

\begin{prop}\label{proposition 1}
Let $ \A = \{0, 1, \cdots, d_1 \} $ be an alphabet; $\A^{\N}$ is endowed with the alternate order. We consider an infinite alternate Lyndon word $(d_i)_{i \geq 1}$ and $ M $ the associated system. $ H_n$ denotes the number of words of length $ n $ of the language of $M$. Then, $ H_0 = 1$ and for all $ n \geq 1 $, we have : 

\begin{equation}  
H_n = \underset{ k = 1 }{\overset{n}{\sum}} (-1)^k(d_{k-1}-d_k)H_{n-k}+1
\label{(5)}
 \end{equation}
with $ d_0 =0$. 
\end{prop}
The proof of the proposition above is deduced from the following lemma:
\begin{lem}\label{lemma 1}
We consider an infinite alternate Lyndon word $(d_i)_{i \geq 1}$ and $ M $ the associated system. $ H_n$ denotes the number of words of length $ n $ in $L_M$.
Let $ A(n)= a_1 a_2 \cdots a_n $ and $ B(n)=b_1 b_2\cdots b_n $ be two words of $ L_M$. Suppose $ A(n) \prec B(n)$. We note $ [A(n), B(n)] $ the set of words $ X(n) = x_1 \cdots x_n \in L_M$ such that 
\begin{equation*}
 A(n) \preceq X(n) \preceq B(n);
\end{equation*}
 $ \Gamma_{A,B}(n)$, the cardinal of $ [A(n), B(n)] $. Then, 
\begin{equation}
\Gamma_{A,B}(n) = \sum\limits_{i=1}^{n} (-1)^i(b_i-a_i)H_{n-i} +1.\label{(100)}
\end{equation}
\end{lem} 

\begin{rem}\label{rem 4}
Let $(d_i)_{i \geq 1}$ be an alternate Lyndon word, $L_M $ the language of the associated Lyndon system $ M $ and $ a\in \{1, \cdots, d_1 \}$. In the meaning of the alternating order, $ ad_1 d_2 \cdots d_{n-1} $ is the largest word with length $ n $ in $L_M $ starting by $ a $ and $ (a-1)0d_1\cdots d_{n-2} $ is the least of length $ n $ starting by $ a-1$. Moreover, $ ad_1 d_2 \cdots d_{n-1} $ and $ (a-1)0d_1\cdots d_{n-2} $ are consecutive in the subset of $ L_M$ of words with length $ n $.
\end{rem}

\begin{preu 4}

 If $ n = 1 $, let $ x_1 \in [A(1), B(1)] $. Thus $ (-1)(a_1 - x_1)\leq 0 $ and $ (-1)(x_1-b_1) \leq 0 $ and then $ a_1 \geq x_1 \geq b_1$, that is $ \Gamma_{A,B}(1) = a_1-b_1 + 1 = (-1)(b_1-a_1)H_0 + 1$. The relation \eqref{(100)} is satisfied for $ n = 1$. Suppose \eqref{(100)} satisfied for $ n \leq k -1 $. 
 
 Suppose $ n=k $.  $ a_1 \cdots a_{k} \prec b_1 \cdots b_{k} $ implies that there exists an integer $ t$ between 0 and $ k $ such that $ a_1 \cdots a_{t} = b_1 \cdots b_{t}$ and $ (-1)^{t+1}(a_{t+1}-b_{t+1})<0$. 
 \begin{itemize}
 \item If $ t\geq 1 $, the cardinal of $ [A(k), B(k)] $ is equal to the number of words between $ a_{t+1} \cdots a_{k} $ and $b_{t+1} \cdots b_{k}$.  
 \begin{equation*}
 \begin{cases}
  a_{t+1} \cdots a_{k} \prec b_{t+1} \cdots b_{k} &\text{ if  $ t $ is even } \\
  b_{t+1} \cdots b_{k} \prec a_{t+1} \cdots a_{k} &\text{ if  $ t $ is odd}
  \end{cases}
 \end{equation*}
 Since \eqref{(100)} is satisfied for $ n \leq k-1$ and the length of $ a_{t+1} \cdots a_{k} $ and $ b_{t+1}  \cdots b_{k} $ is $ k-t \leq k-1 $, we have
\begin{equation*} 
\Gamma_{A,B}(k) = \begin{cases}
                 \sum\limits_{i=1}^{k-t}(-1)^i(b_{t+i}-a_{t+i})H_{k-t-i} +1 &\text{ if $ t $ even}\\
                 \sum\limits_{i=1}^{k-t}(-1)^i(a_{t+i}-b_{t+i})H_{k-t-i} +1 &\text{ if $ t $ odd }
                 \end{cases} 
 \end{equation*}
In the both cases, 
\begin{align*}
\Gamma_{A,B}(k) &= \sum\limits_{i=1}^{k-t}(-1)^{i+t}(b_{t+i}-a_{t+i})H_{k-t-i} +1 \\
                  &= \sum\limits_{i=t}^{k}(-1)^i(b_{i}-a_{i})H_{k-i} +1\\
                  &= \sum\limits_{i=1}^{k}(-1)^i(b_i-a_{i})H_{k-i} +1
\end{align*}
We obtain the last equality thanks to the condition $ a_i = b_i $ for $ 1 \leq i \leq t$.
\item Suppose  $ t = 1$, that is $ a_1 > b_1$. Note that $ a_1d_1 \cdots d_{k-1} \in L_M $ and $ b_1 0 d_1 \cdots d_{k-2} \in L_M$.  
\begin{align*}
a_1\cdots a_{k} & \preceq a_1d_1\cdots d_{k-1} \\ (a_1-1)0d_1 \cdots d_{k-2} &\prec (a_1-1)d_1 \cdots d_{k-1} \\  &\vdots \\ (b_1+1)0d_1 \cdots d_{k-2} &\prec (b_1+1)d_1 \cdots d_{k-1} \\ b_1 0 d_1\cdots d_{k-2} &\preceq b_1 b_2 \cdots b_{k}
\end{align*}
Since $ ad_1\cdots d_{k-1} $ and $ (a-1)0d_1 \cdots d_{k-2} $ are consecutive, it follows that:

\begin{equation*}
[A(k), B(k)]=[A(k), \hspace{0.1cm} a_1d_1\cdots d_{k-1}]\cup [b_10d_1\cdots d_{k-2}, B(k)]\cup 
              \underset{i=b_1+1}{\overset{a_1-1}{\bigcup}}[i0d_1\cdots d_{k-2}, id_1\cdots d_{k-1}] 
\end{equation*}
The intervals $ [i0d_1\cdots d_{k-2}, id_1\cdots d_{k-1}] $ and $ [d_1\cdots d_{k-1}, 0d_1\cdots d_{k-2}] $ have same cardinal 
\begin{equation*} 
H_{k-1}=\sum\limits_{i=1}^{k-1}(-1)^i(d_{i-1}-d_i)H_{k-1-i}+1 
\end{equation*}
(since \eqref{(100)} is supposed satisfied for $ n \leq k-1$). Thus, 
\begin{align*}
\Gamma(k)&= \sum\limits_{i=2}^{k}(-1)^i[(d_{i-1}-a_i)+(b_i-d_{i-2})]H_{k-i}+2+(a_1-b_1-1)H_{k-1}\\
           &= \sum\limits_{i=2}^{k}(-1)^i(b_i-a_i)H_{k-i}+1+(a_1-b_1)H_{k-1} \\
           &= \sum\limits_{i=1}^{k}(-1)^i(b_i-a_i)H_{k-i}+1
\end{align*}
\end{itemize}
\end{preu 4}

\begin{defi}
A language $ L $ over an alphabet $ A $ is factorial  if it contains the sub-words of its words. It is extendable if for all word $ x_1x_2 \cdots x_n $ in $ L $, there exist two letters $ a $ and $ b $ in $ \A $ such that $ ax_1 x_2\cdots x_n b $ belongs to $ L $.
\end{defi}

The language $ L_{M} $ is factorial and extendable. In fact, the extendability to the right  is clear by construction of $ L_M$. For the extendability to the left, remark that for all word $ w \in L_M$, $ 0w \in L_M$.
 That means \textit {$ L_M $ is a language of a dynamical system}. 

We set 
\begin{equation}
h(L_M)=\lim\limits_{n\rightarrow + \infty} \frac{1}{n} \log \sharp (L_{M}\cap \A^n) =\lim\limits_{n\rightarrow + \infty} \frac{1}{n} \log H_n
\label{(14)}
\end{equation}
and $ h(L_{M}) \geq 0$. Since $ L_{M} $ is factorial, according to F. Blanchard and G. Hansel (see \cite{MR858689}), the limit $ h(L_{M}) $ exists. It defines the topological entropy of the system. Note $ \beta $ the positive real number such that 
\begin{equation*}
 \lim\limits_{n\rightarrow+\infty} \frac{1}{n} \log \sharp (L_{M}\cap \A^n) = \log \beta. 
\end{equation*}

\begin{rem}\label{rem 5}
The real $\beta$ defined above is between $ d_1 $ and $ d_1 + 1$. Indeed, the full-shift $ \{0, 1, \cdots, d_1 \}^{\N} $ which, endowed with the shift $ \sigma $ has $\log (d_1 + 1)$ as entropy, contains the alternate Lyndon system associated to $\beta$. Moreover, the alternate Lyndon system (Lyndon system of alternating Lyndon word) associated  to $ \beta $ contains strictly the full-shift $ \{0, 1, \cdots, d_1 - 1\}^{\N}$ which, endowed with the shift is of entropy $ \log d_1$. 
\end{rem}

\begin{prop}\label{proposition 2}
Let $ (a_i)_{i \geq 1} $ and $ (d_i)_{i \geq 1} $ be two infinite alternate Lyndon words over an alphabet $ \A $ such that :
\begin{equation*}
d_1d_2d_3\cdots \prec a_1 a_2 a_3 \cdots.
\end{equation*}

Note $ M^{'} $ the alternate Lyndon system associated to $ (a_i)_{i \geq 1} $ and $ M $ that associated to $(d_i)_{i \geq 1} $. Then, $  M $ contains $ M^{'} $ and so, $ h(L_{M^{'}}) \leq h(L_{M})$. 
\end{prop}

About the entropy of the system, two situations can occur: $M$ has an entropy equal to zero or this one is strictly greater than 0. In the paragraphs above, we justify the existence of the entropy of $ M $. The following subsection determines the systems of type $ M $ with non-zero entropy. In other words, we give an answer at the question : for what alternate Lyndon word $ (d_i)_{i \geq 1}$ the real $\beta$ is strictly greater than 1?

We define on the alphabet  $ \left\lbrace 0, 1 \right\rbrace  $ the morphism $\phi $ by : $ \phi (0) = 1 $ and $ \phi (1) = 100 $. We set $ \phi^{\infty} (1) = \lim\limits_{n \rightarrow +\infty} \phi^n (1) $.

\begin{equation}
\phi^{\infty}(1) = 1001110010010011100111001110010010011100100\cdots .
\label{(15)}
\end{equation}
Replacing 1 by 2 and 0 by 1, we have :

\begin{equation*}
s=211 2 2 211 211 21122 21122 21122211211 21122211211 \cdots .
\end{equation*}
Put in this form, we recognize in $ w = 1s $ the sequence A026465 of the encyclopedia of Sloane. The elements of $ w $ count the number of consecutive identical symbols in the Thue-Morse sequence $ t $ defined by : 

\begin{equation*}
t=(t_i)_{i \geq 0} = 0110100110010110 \cdots .
\end{equation*}
It is obtained by the induction $ t_0= 0$, $ t_{2n} = t_n $ and $ t_{2n+1} = 1-t_n$, or by successive iteration of the morphism over $ \{0, 1\} $ defined by $ 0 \mapsto 01 $ and $ 1 \mapsto 10 $. 

We set $ u_n = \phi^{n}(1) $, $ v_0 = 00 $ and for all $ n> 0 $,  $ v_n = u_{n-1}u_{n-1} $. So, 
\begin{equation*}
\begin{aligned}
u_n &= u_{n-1} v_{n-1} \\
    &= u_{n-1} u_{n-2} u_{n-2}.
\end{aligned}
\end{equation*}

\begin{theo}\label{theorem 1}
Let $ (d_i)_{i \geq1}$ be a Lyndon word over a finite alphabet $\A=\{0, 1, \cdots, d_1 \} $. We note $ M $ the associated Lyndon system. Then,
\begin{itemize}
\item[(a)] $ M $ has a positive entropy if only if
\begin{equation*}
(d_i)_{i \geq 1} \prec \phi^{\infty}(1).
\end{equation*}
\item[(b)] $ M $ has zero entropy if and only if 
\begin{equation*}
(d_i)_{i\geq 1}\in \{\overline{0},\phi(1)\}\cup \{\overline{\phi^n(1)} : n\in \N \}.
\end{equation*}
\end{itemize}
\end{theo}

\subsection[System with positive entropy]{Alternate Lyndon system with non zero entropy}

\begin{prop}\label{proposition 3}

Let given a Lyndon system for the alternate order $ M((d_i)_{i \geq 1}) = M $, $ H_n $ the number of words of length $ n $ of the language. We suppose its entropy equals $ \log \beta > 0$; then $ \beta $ is the largest real solution of 
\begin{equation}
 -\frac{x}{x+1} = \sum\limits_{n \geq 1}\frac{d_n}{(-x)^n}.
 \label{16}
\end{equation}

\end{prop}

\begin{preu}
Since $ \lim \limits_{n \rightarrow +\infty} H_n^{1/n} $ equal to $ \beta $ and $ \beta> 1$, $ \frac{1}{\beta} $ is the radius of convergence of the power series
$ \sum \limits_{n \geq 0} H_n z^n $. In the open disk of center 0 and radius 
 $ 1/\beta $, from \eqref{(5)}, we have

\begin{align*}
 \sum\limits_{n \geq 0} H_n z^n & = \sum\limits_{n \geq 0} \left( \sum\limits_{k = 1}^n(-1)^k(d_{k-1}-d_k)H_{n-k} \right) z^{n} + \sum\limits_{n \geq 0} z^n\\
                                   & = \left( \sum\limits_{n \geq 0} H_n z^n \right) \left( \sum\limits_{n \geq 1} (d_{n-1}-d_n)(-z)^n \right) + \sum\limits_{n \geq 0} z^n.
\end{align*}
Then

\begin{equation}
\sum\limits_{n\geq 0} H_n z^n = \dfrac{1}{(1-z)(1-(1+z)\sum\limits_{n \geq 1}d_n(-z)^n)}.
 \label{(18)}
\end{equation}
It follows from Pringsheim’s theorem (see theorem IV.6 in \cite{MR2483235}) that $\frac{1}{\beta}$ is a singularity of $ \sum\limits_{n \geq 0} H_nz^n$ and thus a root of $1-(1+z)\sum\limits_{n \geq 1} d_n(-z)^n$.
\end{preu}

\begin{lem}\label{lemma 2}
For all $ n $ in $ \N $, $ u_n = \phi^n(1)  $ is a word over $ \{0, 1\} $ with odd length. 
\end{lem}
For the proof of the previous lemma, see theorem 2.6 of \cite{ MR2974214}

\begin{prop}\label{proposition 4}
For all $ n \in \N $, $ u_n \overline{ v_n } $ is a Lyndon word with respect to the alternate order.
\end{prop}
The previous proposition is a consequence of theorem 2.5 of \cite{MR2974214}. Furthermore, for all $ k \in \N$, $ \vert u_k \vert = 2 \vert u_{k-1} \vert - (-1)^k $ and $ \vert u_n \vert $ and $ \vert v_n \vert $ are consecutive integers.

Let $ a_1 a_2 \cdots a_n $ be a word of a finite alphabet $\A $. We set 
\begin{equation*}
(1a_1a_2 \cdots a_n) (z) = 1+\sum\limits_{k=1}^na_k(-z)^k.
\end{equation*}
We define on $ \{1,2\}^{*}$ the substitution $ \psi $ such that $\psi(1) = 2$ and $ \psi(2) = 211 $. Note $ \omega_n = \psi^n(2) $. From the proof of lemma 5 of \cite{MR2292537},
\begin{equation}
(1w_n)(z) = \prod\limits_{k = -1}^{n-1} (1-z^{\vert u_k \vert}) - z^{\vert u_n \vert}
\label{M}
\end{equation}
with $ u_{-1} = 0 $. The relation \eqref{M} can be obtained by induction on $ n $.

\begin{prop}\label{proposition 5}

\begin{itemize}

 \item[(1)] For all $ n \in \N $, the Lyndon system associated to the alternate Lyndon word $ u_n \overline{v_n} $ has non zero entropy.

 \item[(2)] The Lyndon system $ M(\phi^{\infty}(1))$ associated to the alternate Lyndon word $ \phi^{\infty}(1) $ is of entropy zero.
 \end{itemize}
\end{prop}

\begin{preu}

\begin{itemize}
\item[(1)] $ M(u_n\overline{v_n} )$ has non zero entropy :
\end{itemize}
 The power series $ 1+ \sum\limits_{n \geq 1} (-1)^n(d_n+1)z^n $ (where $ d_1 d_2 \cdots = u_n \overline{v_n} $) converges in the ball of radius 1 and center 0. The entropy of the system is the logarithm of the inverse of the smallest zero (in modulus) of this power series (see proposition \ref{proposition 3}). Since $ (d_1+1)(d_2+1)\cdots (d_n+1) \cdots = \omega_n\overline{\omega_{n-1}}$ and from \eqref{M}, 

\begin{eqnarray*}
 1+ \sum\limits_{n \geq 1}(d_n+1)(-z)^n= \frac{1}{1+z^{\vert u_{n-1} \vert}}\left( \prod\limits_{k=-1}^{n-2} (1-z^{\vert u_k \vert}) \right) (1-z^{2\vert u_{n-1} \vert} - z^{\vert u_n \vert}) 
 \end{eqnarray*}
 where $ u_{-1} = 0 $. Then, the smallest zero $ \frac{1}{\gamma_n}$ satisfies :
 \begin{equation*}
  1 = \frac{1}{\gamma_n^{\vert u_n \vert }} + \frac{1}{ \gamma_n^{\vert v_n \vert }},
 \end{equation*}
That is, 
\begin{equation*}
\gamma_n^{l_n } = \gamma_n  + 1 \text{ with $ l_n = \max ( \vert u_n \vert, \vert v_n \vert ) $ },
\end{equation*}
This implies that $ \log \gamma_n > 0 $.

\begin{itemize} 
 \item[(2)] $ M(\phi^{\infty}(1)) $ has zero entropy : 
\end{itemize} 
The sequence $ (u_n \overline{v_n})_{n \geq 1} $ is increasing in the sense of the alternate order.
 \begin{equation*}
 u_n\overline{v_n} \prec u_{n+1}\overline{v_{n+1}} \prec \cdots \prec \phi^{\infty}(1).
 \end{equation*}
 Indeed, $ u_n \overline{v_n} = u_{n+1}\overline{v_n} $ and $ \overline{v_{n+1}} \prec \overline{v_{n}} $ since $ v_{n+1} = u_nu_n\prec v_nv_n $ and add $ u_{n+1}$ at left of $ u_n u_n $ and $ v_n v_n $. We note $ M_{\gamma_n}$ the system associated to $ u_n \overline {v_n}$, with $ \log \gamma_n $ the entropy. We have a sequence of Lyndon systems $ (M_{\gamma_n})_{n \geq 1}$ such that

\begin{equation*}
 M(\phi^{\infty})\varsubsetneq \cdots \varsubsetneq  M_{\gamma_{n+1}} \varsubsetneq M_{\gamma_n}\varsubsetneq \cdots\varsubsetneq  M_{\gamma_0}  
 \text{ and $ \gamma_{n+1} < \gamma_n < \cdots < \gamma_0$ }, 
\end{equation*}
and
\begin{equation*}
  \gamma_n^{l_n} = \gamma_n + 1
\end{equation*}
where $ l_n = \max ( \vert u_n \vert, \vert v_n \vert ) $. 

The length $ \vert u_n \vert $ and $ \vert v_n \vert $ tend to infinity with $ n $. It is the same for $ l_n $. It follows that $ \gamma_n $ tends to 1. Then, $ M(\phi^{\infty}) $ has zero entropy. 
\end{preu}

As consequence of this proposition, we have the following corollary :

\begin{coro}\label{coro 1}
A system $ M $ associated to the alternate Lyndon word $(d_i)_{i \geq 1} $ has a  strictly positive entropy if only if
 
 \begin{equation}
 (d_i)_{i \geq 1} \prec \phi^{\infty}(1).
 \label{(21)}
 \end{equation}
\end{coro}

\begin{preu}

 (a) Suppose $M$ with non-zero entropy. In this case $ (d_i)_{i \geq 1} \neq \lim \limits_{n \rightarrow +\infty } \phi^n(1) $.
From the proposition \ref{proposition 2}, $ (d_i)_{i \geq 1} \prec \lim \limits_{n \rightarrow +\infty} \phi^n(1)$. Otherwise, $ M$ has zero entropy.

 (b) Suppose $ (d_i)_{i \geq 1} \prec \lim \limits_{n \rightarrow +\infty} \phi^n(1) $. Since $ \lim \limits_{n \rightarrow +\infty} \phi^n(1) = \lim \limits_{k \rightarrow +\infty} \phi^k(1) \overline{\phi^{k-1}(1)}$, there exists an integer $m$ such that $ d_1 d_2 d_3 \cdots \prec \phi^m(1) \overline{\phi^{m-1}(1)}$. So, $ M $ contains the system associate to the alternate Lyndon word $ u_m\overline{v_m} $ which has $ \log \gamma_m $ as non zero entropy. Then, $ M $ has a positive entropy.
\end{preu}

\subsection[System with entropy zero]{Alternate Lyndon systems with zero entropy}

The alternate Lyndon word $ \overline{10} $ is attached to 2. The associated system has an entropy equal to $\log 2 > 0$. It follows that all alternate Lyndon system with zero entropy is attached to an alternate Lyndon word over the alphabet $ \{0, 1\}$.

\begin{rem}\label{rem 6}
\begin{enumerate}
\item $ 1111 \cdots = \overline{1} $ is the alternate Lyndon word (different to  $ \overline{0} $) generating the smallest alternate Lyndon system. That is,
 
\begin{equation*}
d_1 d_2 d_3 \cdots \preceq 111111\cdots, \text{ and $ (d_i)_{i \geq 1} \neq \overline{0} $}.
\end{equation*}

\item Let $ u $ be a word with odd length over a finite alphabet. It is easy to show that if $ (d_i)_{i \geq 1} $ begins by the sequence $ u u $ then, $ (d_i)_{i \geq 1} = \overline{u}.$
\end{enumerate}
\end{rem}

\begin{preu 3}
We obtain the first assertion of the theorem \ref{theorem 1} thanks to the proposition \ref{proposition 5} and corollary \ref{coro 1}. Let prove (b).

Let $ M \neq M(\phi^{\infty}(1)) $ be a Lyndon system attached to an alternate Lyndon word $ (d_i)_{i \geq 1} $. $ M $ has zero entropy if only if

\begin{equation*}
 \phi^{\infty}(1) \prec (d_i)_{i \geq 1}.
\end{equation*}
There exists an integer $ n $ such that : 
  \begin{equation*}  
  u_{n} = d_1d_2 \cdots d_{\vert u_{n} \vert} \text{ and $ u_n \prec d_{\vert u_n \vert + 1} \cdots d_{\vert u_{n} \vert + \vert v_n \vert} \prec v_n $ }.
  \end{equation*}
  Remark that : 
  \begin{equation*} 
  u_{2k} = u_{2k-1} u_{2k-2} \cdots u_1 u_0 u_0 \text{ and } v_{2k} = u_{2k-1} \cdots u_1u_0 00.
  \end{equation*}
Similarly,
\begin{equation*} 
  u_{2k+1} = u_{2k} u_{2k-1} \cdots u_1u_0 00 \text{ and } v_{2k+1} = u_{2k} u_{2k-1} \cdots u_1 u_0 u_0.
\end{equation*} 
In the both cases, $ u_n $ and $ v_n $ begin by : $ u_{n-1} u_{n-2} \cdots u_1u_0 $. Set
\begin{equation*}  
W=d_{\vert u_n \vert + 1} \cdots d_{\vert u_{n} \vert + \vert v_n \vert}. 
\end{equation*}
Then, $ W $ begins by 
 $ u_{n-1} u_{n-2} \cdots u_1u_0 $.
\begin{itemize} 
\item if $ n $ is even, $ (-1)^{\vert u_n \vert} w_{\vert u_n \vert }< 0$ where $ w_{\vert u_n\vert} $ denotes the $ \vert u_n \vert$-th letters of $ W $. In fact, at the index $ \vert u_n \vert $ of $ v_n $, there is 0. It follows that $ w_{\vert u_n \vert} = 1 =u_0 $. Hence, $ W $ begins by 
\begin{equation*}
u_{n-1}u_{n-2} \cdots u_1 u_0 u_0=u_n .
\end{equation*}
So, $ (d_i)_{i \geq 1}$ begins by $ u_n u_n $. Considering the previous remark, $ (d_i)_{i \geq 1} = \overline{u_n} $.
\item Similarly, if $ n $ is odd, $w_{\vert v_n \vert } = 0 $. Then, $ W=u_{n-1} \cdots u_1u_0 0$. But $ u_0 $ is always followed by 00. Thus, 

\begin{equation*}
 d_{\vert u_n \vert + 1} \cdots d_{2\vert u_n \vert} = u_n . 
\end{equation*}
\end{itemize}
From the previous remark, we conclude that $ (d_i)_{i \geq 1}= \overline{u_n} $. 
\end{preu 3}

An immediate consequence of the theorem \ref{theorem 1} is the lemma 4 of \cite{MR2204744} : $ \phi^{\infty}(1) $ is the greatest (with respect to the alternate order) non periodic alternate Lyndon word. That is :
\begin{equation*}
\forall (d_i)_{i \geq 1} : (d_i)_{i \geq 1} \prec (d_{i+k})_{i \geq 1}, \text{ $ \forall k \in \N $},
\end{equation*}
then,
\begin{equation*}
(d_i)_{i \geq 1} \preceq \phi^{\infty}(1).
\end{equation*}
In other words, $M(\phi^{\infty}(1))$ is the smallest system attached to a non periodic alternate Lynodn word. Moreover, $ M(\phi^{\infty}(1)) $ constitutes the greatest system with zero entropy.

In the following, we are interested only non-zero entropy alternate Lyndon systems. That is, we will consider the alternate Lyndon words $(d_i)_{i \geq 1} $ such that 
\begin{equation*}
 (d_i)_{i \geq 1} \prec \lim \limits_{n \rightarrow + \infty} \phi^n(1) .
 \end{equation*}

\section{An increasing map}

Let $ M $ be a Lyndon system attached to an alternating Lyndon word, $ \log \beta $ the entropy and $ L_M $ the language of its language. We define the map $f_{\beta} $ from $ L_M $ to $ \R $ by 
\begin{equation*}
f_{\beta} (x_1x_2 \cdots x_n) = \sum\limits_{n \geq 1}^n\frac{x_k}{(-\beta)^k}.
\end{equation*}
The map $ f_{\beta }$ can be defined on infinite words of $ M $. For $ (x_i)_{i \geq 1} $ in $ M $, 
\begin{equation*}
f_{\beta} ((x_i)_{i \geq 1}) = \sum\limits_{n \geq 1}\frac{x_n}{(-\beta)^n}.
\end{equation*}
\begin{lem}\label{lemma 3}

Let $ (a_i)_{i \geq 1} $ and $ (d_i)_{i \geq 1} $ be two infinite alternate Lyndon words over a finite alphabet $ \A $ such that : 

\begin{equation*}
d_1d_2d_3\cdots \prec a_1 a_2 a_3 \cdots . 
\end{equation*}
Let  $ \log \beta $ be the entropy of the system associated to $(d_i)_{i \geq 1}$ that we note $ M $. Then:

\begin{equation}
\sum\limits_{ n \geq 1} \frac{a_n}{(-\beta)^n} \geq -\dfrac{\beta}{\beta + 1}.
\label{(22)}
\end{equation}

\end{lem}

\begin{preu}
The function $ x \mapsto - \frac{x}{x+1} $ on $ [1, +\infty ) $ decreases until $ - 1 $. From proposition \ref{proposition 3} , $ \beta$ is the largest real number such that $ \sum\limits_{ n \geq 1} \frac{d_n}{(-\beta)^n} = -\dfrac{\beta}{\beta + 1} $. Let $ M^{'} $ be the system associated to $ (a_i)_{i \geq 1}$ and $ \log \alpha $ its entropy. $ \alpha $ is the largest real number satisfying $ \sum\limits_{n\geq 1}\frac{a_n}{(-\alpha)^n} = -\dfrac{\alpha}{\alpha + 1} $.  

 The graph of $ x \mapsto \sum\limits_{ n \geq 1} \frac{a_n}{(-x)^n} $ intersects the one of $ x \mapsto - \frac{x}{x+1} $ in $ \alpha $. Moreover, 
 \begin{equation*} 
 \sum \limits_{ n \geq 1} \frac{a_n}{(-x)^n}  \underset{ x \rightarrow + \infty }{\longrightarrow} 0  
 \end{equation*}
and 
\begin{equation*}
-\frac{x}{x+1} \underset{x\rightarrow +\infty}{\longrightarrow} -1 .
\end{equation*}
 Then, after $ \alpha $, the graph of $\sum\limits_{ n \geq 1} \frac{a_n}{(-x)^n} $ is above the one of $ - \frac{x}{x+1} $. 
 
$d_1 d_2 \cdots \prec a_1 a_2 \cdots $ implies that $ M^{'}$ is included in $ M $ and then, the entropy of $ M^{'}$ is less than the one of $ M $, that is $ \alpha \leq \beta $. Hence the result. 
  \end{preu}
  
 One of the consequences of the previous lemma is given in following corollary:
 \begin{coro}\label{coro 2}
 Let $(d_i)_{i \geq 1} $ be an alternating Lyndon word, $ M $ the associated system, $ L_M $ the language of $ M $ and $\log \beta $ its entropy. Suppose $ d_1 d_2 \cdots d_{i-1} j \in L_M$ such that $ (-1)^i(d_i-j)< 0 $. Then, 
 \begin{equation*}
 f_{\beta}(d_1 d_2 \cdots d_{i-1} j) \geq -\dfrac{\beta}{\beta+1}.
 \end{equation*}
 \end{coro}
\begin{preu}

Remark that if $ j\neq d_1 $, $d_1 d_2 \cdots d_{i-1} j \overline{0} $ is an alternating Lyndon Word satisfying the condition $ d_1d_2\cdots \prec d_1 d_2 \cdots d_{i-1} j \overline{0} $ (since $ d_1 d_2 \cdots d_{i-1} j \in L_M$) and we apply the previous lemma. 

If $ j = d_1 $, that is $ i $ is even, we consider the alternating Lyndon word $ \overline{d_1\cdots d_{i-1}d_1} $ and $ d_1d_2 \cdots \preceq  \overline{d_1\cdots d_{i-1}} $. Hence, 
\begin{equation*} 
f_{\beta}(\overline{d_1\cdots d_{i-1}d_1} ) = f_{\beta}(d_1\cdots d_{i-1}d_1) \frac{\beta^{i}}{\beta^{i}-1} \geq -\frac{\beta}{\beta+1} 
\end{equation*} 
The integer $i$ is even. It follows that $ f_{\beta}(d_1\cdots d_{i-1}d_1) \geq -\frac{\beta}{\beta+1} $ (this inequality holds because $\frac{\beta^{i}}{\beta^{i}-1} > 1 $).
\end{preu}

\begin{rem}\label{rem 7}
If $ d_1 \cdots d_{k}(d_{k+1}-(-1)^k) \in L_M $ and $ d_1 \cdots d_{k-1}(d_{k}+(-1)^k) \in L_M $, from the previous corollary, 
\begin{equation*}
f_{\beta}(d_1 \cdots d_{k+1})-f_{\beta}(0d_1 \cdots d_{k}) \geq -1.
\end{equation*}
And then, 
\begin{equation}
f_{\beta}(ad_1 \cdots d_{k+1})<f_{\beta}((a-1)0d_1 \cdots d_{k}) \label{O}. 
\end{equation}
The word $ d_1\cdots d_{k}(d_{k+1}-(-1)^k) \not\in L_M $ if $ k $ is even and $ d_{k+1} = 0$ or if there exists and even integer $ i $ between $ 1 $ and $ k+1 $ such that
\begin{equation}
d_1 \cdots d_{k+1} = d_1 \cdots d_{i+1}d_1 \cdots d_{k-i} \label{P}.
\end{equation} 
If $ k $ even and $ d_{k+1} = 0$, 
\begin{equation}
f_{\beta}(ad_1\cdots d_{k}0)-f_{\beta}((a-1)0d_1 \cdots d_{k}) \leq f_{\beta}(ad_1\cdots d_{k})-f_{\beta}((a-1)0d_1 \cdots d_{k-1}) \label{X}.
\end{equation}
\end{rem}

\begin{prop}\label{proposition 7}
Let $ x_1 x_2 \cdots x_n $ and $ y_1 y_2 \cdots y_n $ be two finite words of the language $M$. Then,  
 \begin{equation} 
 x_1x_2 \cdots x_n \prec y_1 y_2 \cdots y_n  \Rightarrow f_{\beta}( x_1x_2 \cdots x_n ) < f_{\beta}(y_1 y_2 \cdots y_n). \label{222} 
 \end{equation}
\end{prop}
\begin{preu}

 For all finite word $ a_1 a_2 \cdots a_n\in L_M $, we have : 

\begin{equation*}
(d_1, d_2, \cdots, d_n) \preceq (a_1, a_2, \cdots, a_n) \preceq (0, d_1, d_2, \cdots, d_{n-1}).
\end{equation*}

\begin{equation*}
 (x_1) \prec (y_1) \Rightarrow 
                -\frac{x_1}{\beta} < -\frac{y_1}{\beta}. 
\end{equation*}

The property is true for $ n = 1 $. Suppose it true for $n \leq k $.

For $n= k+1 $, let $( x_1, x_2, \cdots, x_{k+1}) \prec (y_1, y_2, \cdots, y_{k+1}) $. So, there exists $ j $ such that $ x_i = y_i $ for all $ i < j $ and $ (-1)^j(x_j-y_j)< 0$.

If $ j \geq 2 $, the length of $ (x_j, x_{j+1}, \cdots, x_{k+1}) $ and $ (y_j, y_{j+1}, \cdots, y_{k+1}) $ is less than or equal to $k$ and we obtain immediately \eqref{222}. In fact, Thanks to \eqref{111},
\begin{equation*}
\begin{cases} 
f_{\beta}(x_j \cdots x_{k+1}) \prec f_{\beta}(y_j\cdots y_{k+1}) &\text{ if $ j $ is odd }\\
f_{\beta}(y_j \cdots y_{k+1}) \prec f_{\beta}(x_j\cdots x_{k+1}) &\text{ if $ j $ is even}
\end{cases}
\end{equation*}
In the both cases, 
\begin{equation*}
\frac{1}{(-\beta)^{j-1}}f_{\beta}(x_j \cdots x_{k+1}) < \frac{1}{(-\beta)^{j-1}} f_{\beta}(y_j \cdots y_{k+1}). 
\end{equation*}
And then,
\begin{equation*}
 f_{\beta}(x_1 \cdots x_{k+1}) < f_{\beta}(y_1 \cdots y_{k+1}) 
\end{equation*}
since $ x_t=y_t$ for $1\leq t \leq i-1$.

Note that for $ j=1 $, $ x_1 \geq y_1 + 1 $. So, the language $ L_M $ of $ M $ contains $ x_1 d_1 d_2 \cdots d_k $ and $ y_1 0 d_1 d_2 \cdots d_{k-1} $. 
The maximum (with respect to alternate order) of length $k+1$ beginning by $x_1$ is $ x_1 d_1 d_2 \cdots d_k $ and the minimum (with respect to the alternate order) of length $k+1$ beginning by $y_1$ is $y_10d_1\cdots d_{k-1}$. Set $ X^{(k+1)}=x_1 x_2 \cdots x_{k+1} $ and $ Y^{(k+1)} = y_1 \cdots y_{k+1} $. Hence, 
\begin{equation*}
 X^{(k+1)} \preceq x_1 d_1 \cdots d_k \prec (x_1-1)0d_1 \cdots d_{k-1} \preceq y_1 0 d_1 \cdots d_{k-1} \preceq Y^{(k+1)}.
\end{equation*}

From the previous case, 
$f_{\beta}(x_1 x_2 \cdots x_{k+1})$ and $f_{\beta}(y_1 0 d_1 \cdots d_{k-1}) $ are less than $ f_{\beta}(x_1 d_1 \cdots d_k)$ and 
 $ f_{\beta}(y_1 \cdots y_{k+1})$ respectively. We set  
\begin{equation*}
 \Delta(k) = f_{\beta}(x_1 d_1 \cdots d_k)- f_{\beta}((x_1-1)0d_1 \cdots d_{k-1}).
 \end{equation*}

If there exists $ i $ odd ($ 1 \leq i \leq k $) such that \eqref{P} is hold.
\begin{align*}
\Delta = &\Delta(i) - \frac{1}{\beta^{i+1}}\left(f_{\beta}(d_1\cdots d_{k-i})-f_{\beta}(d_{i}d_1\cdots d_{i-1}) \right) \\
         &< 0
\end{align*}
Using  \eqref{O} and \eqref{X}, we complete the proof that $ \Delta < 0$. 

Since
$ f_{\beta}((x_1 - 1)0d_1\cdots d_{k-1}) \leq f_{\beta}(y_10d_1 \cdots d_{k-1})$, 
we conclude that $ f_{\beta}$ is an increasing map on finite words.
 \end{preu}
\vspace{1cm}

\begin{theo}\label{theorem 2}
Let $(d_i)_{i \geq 1} $ an infinite alternate Lyndon word, $ M = M((d_i)_{i \geq 1}) $ the associated one-side dynamical system. We set $ \log \beta $ its entropy and $ f_\beta$ the map from $ M $ defined by 
\begin{equation*}
f_\beta((a_i)_{i \geq 1}) = \sum\limits_{n \geq 1}\dfrac{x_n}{(-\beta)^n}.
\end{equation*} 
Then,
 \begin{itemize}
  
  \item $ f_{\beta} $ is increasing on $ M $. 
  
  \item $ f_{\beta} $ is continuous on $ M $.
  
  \item The image by $ f_{\beta} $ of $ M $ is the interval $\overline{I}_{\beta} = [ - \frac{\beta}{\beta + 1}, \frac{1}{\beta + 1} ] $.
 \end{itemize}

\end{theo}

\begin{preu}
\begin{itemize}
\item The first assertion of this theorem results of the previous proposition. The  growth on the infinite words is obtained by taking the limit on $ n $ in \eqref{222} (the word are over a finite alphabet) and the continuity of $ f_{\beta}$ is obvious. 
\end{itemize}

\begin{itemize}
\item \textbf{Image of $ M $ by $ f_{\beta}$}.
\end{itemize}
Let show that the image by $ f_\beta$ of $M$ is the closed interval $ \overline{I}_{\beta} $ with bounds $ -\dfrac{\beta}{\beta + 1} $ and $ \dfrac{1}{\beta + 1}$.

Since $ f_{\beta} $ is  an increasing map, $ \overline{I}_{\beta} $ contains $ f_{\beta}(M)$. If $ \overline{I}_{\beta} \neq f_{\beta}(M) $, $ \overline{I}_{\beta} \setminus f_{\beta}(M) $ contains an open interval $ (r, t) $ such that there exists two sequences $ (x_i)_{i \geq 1} $ and $ (y_i)_{ i \geq 1} $ in $ M $ satisfying :
\begin{equation*}
 f_{\beta}(x_1x_2 \cdots) = r \text{ and } f_{\beta}(y_1 y_2 \cdots)=t.
\end{equation*}
\begin{equation*}
  ( r, t ) \subset \overline{I}_{\beta} \setminus f_{\beta} (M),   
\end{equation*}
with $ r $ and $t$ consecutive in $ f_{\beta}(M) $.
Considering the properties of the alternate order, $ x_1x_2 \cdots \prec y_1 y_2 \cdots $. The reals $ r $ and $ t $ are different, the equality between the both sequences can not occur. Then, there exists an integer $ n $ such that $ x_i = y_i $ for $ i<n $ and $ (-1)^n(x_n - y_n)< 0 $.

 If $ (-1)^n(x_n - y_n) \leq -2 $, there are as much elements in $ M $ as in the interval of words $ x_1 x_2 \cdots $ and $ y_1 y_2 \cdots $. Indeed, all concatenations of $ x_1x_2 \cdots x_{n-1}(x_n+(-1)^n) $ and a word $ u \in M $ is such that :

\begin{equation*}
 x_1x_2 \cdots x_n x_{n+1} \cdots \prec x_1 x_2 \cdots x_{n-1}(x_n+(-1)^n)u \prec y_1 y_2 \cdots .
\end{equation*}
Since $ u $ sweeps across $ M $, the image by $f_{\beta}$ of the set of these words is $ \gamma + \frac{1}{(-\beta)^n} f_{\beta}(M) $, with 
\begin{equation*}
\gamma = f_{\beta}(x_1x_2 \cdots x_{n-1}(x_n+(-1)^n).
\end{equation*}
This contradicts the fact that $ t $ and $ r $ are consecutive in $ f_{\beta}(M) $. 

Now, consider $ (-1)^n(x_n - y_n) = -1 $. So, 

\begin{equation*}
 r-t = \frac{1}{(-\beta)^n}\left( -(-1)^n+f_{\beta}(x_{n+1}x_{k+2} \cdots) - f_{\beta}(y_{n+1}y_{k+2} \cdots)\right).
\end{equation*}
Then, $ r \neq t $ means :

\begin{equation*}
 \begin{cases}
  f_{\beta}(x_{n+1}x_{n+2}\cdots ) - f_{\beta}(y_{n+1}y_{n+2} \cdots) \neq 1  & \text{ if $ n $ is even }, \\
  f_{\beta}(x_{n+1}x_{n+2}\cdots ) - f_{\beta}(y_{n+1}y_{n+2} \cdots) \neq -1 & \text{ if $ n $ is odd }.
 \end{cases}
\end{equation*}
That is : 

\begin{equation*}
 \begin{cases}
 (f_{\beta}(x_{n+1}x_{n+2} \cdots), f_{\beta}(y_{n+1}y_{n+2}\cdots)) \neq ( \dfrac{1}{\beta + 1}, -\dfrac{\beta}{\beta + 1} ) & \text{ if $ n $ is even }, \\
  ( f_{\beta}(x_{n+1}x_{n+2} \cdots), f_{\beta}(y_{n+1}y_{n+2}\cdots)) \neq ( -\dfrac{\beta}{\beta + 1}, \dfrac{1}{\beta + 1} ) & \text{ if $ n $ is odd}.
 \end{cases}
\end{equation*}

The subshift $ M $ is invariant by the shift. Then, $ x_{n+1}x_{n+2} \cdots $ and $ y_{n+1} y_{n+2} \cdots $ belong to $ M $. So, all word $ v $ beginning by 
$ x_1 x_2 \cdots x_n $ is such that :

\begin{equation*}
\begin{cases}
  x_1x_2 \cdots x_n 0d_1 d_2 \cdots \preceq v \preceq x_1 x_2 \cdots x_n d_1 d_2 \cdots &\text{ if $ n $ is odd }, \\
 x_1x_2 \cdots x_n d_1 d_2 \cdots \preceq v \preceq x_1 x_2 \cdots x_n 0d_1 d_2 \cdots &\text{ if $ n $ is even }.
 \end{cases}
\end{equation*}

Let $ (d^{*}_i)_{i \geq 1} $ be the maximum, in sense of the alternate order, in $ M $ such that :  
\begin{equation*}
f_{\beta}(d^{*}_1 d^{*}_2 d^{*}_3 d^{*}_4 \cdots) = -\dfrac{\beta}{\beta + 1}.
\end{equation*}

We focus on the case $ n $ even, the odd case dealing similarly. If 
\begin{equation*}
 f_{\beta}(x_1 x_2 \cdots) \neq \dfrac{1}{\beta+1}
\end{equation*}
then,
\begin{equation*}
 x_1 x_2 \cdots \prec x_1 x_2 \cdots x_n 0d^{*}_1 d^{*}_2 \cdots \prec x_1x_2 \cdots x_{n-1}(x_n+1)d_1d_2 \cdots \preceq y_1y_2 \cdots .
\end{equation*}
There exists an integer $m$ such that 
\begin{equation*} 
x_{n+1} \cdots x_{n+m-1} = 0d^{*}_1 d^{*}_2\cdots d^{*}_{m-2} \text{ and }  (-1)^m(x_{m+n}-d^{*}_{m-1})< 0.
\end{equation*}
 All words beginning by $ x_1 x_2 \cdots x_n 0 d^{*}_1 d^{*}_2 \cdots d^{*}_{m-1} $ are greater than $ x_1 x_2 x_3 \cdots $ in the sense of the alternate order. Then, it suffices to consider a suitable integer $ i $ greater than $ m $ and consider the words $ x_1 \cdots x_n 0d^{*}_1 \cdots d^{*}_{i-1}j u $ with $ (-1)^{i}(d^{*}_i - j)< 0$ and $ u $ sweeping across $ M $.
The image by $ f_\beta$ of these words contains as much elements as $f_{\beta}(M)$. Then, an infinite numbers of them are less than $ f_{\beta}(y_1 y_2 y_3 \cdots) $. This contradicts the fact that $ t $ and $ r $ are consecutive in $ f_{\beta}(M) $.
\end{preu}

\vspace{1cm}
According to Proposition \ref{proposition 7} and the fact that the map $ f_{\beta} $ is onto, it follows that:
 \begin{prop}\label{proposition 8}
Let $ (d_n)_{n \geq 1} $ be an alternate Lyndon word associated to $\beta > 1$, the images by $ f_{\beta}$ of two cylinders $ _0\left[   x_1 \cdots x_n \right] $ and $_0\left[ y_1 \cdots y_n \right] $ of same length are two disjoint intervals or contiguous of intersection reduced to one point. If these images have in common one point $ z $ and if  
\begin{equation*}
x_1 \cdots x_n \prec y_1 y_2 \cdots y_n \text{ and } x_1 x_2 \cdots x_{n-1} = y_1 y_2 \cdots y_{n-1}.
\end{equation*}
then,
\begin{equation*}
y_1 y_2 \cdots y_n = x_1x_2 \cdots x_{n-1} (x_n + (-1)^n)
\end{equation*}
and
\begin{equation*}
 z = f_{\beta}(x_1 x_2 \cdots x_n d_1 d_2 \cdots) = f_{\beta}(x_1 \cdots x_{n-1}(x_n + (-1)^n)0d_1 d_2 \cdots).
\end{equation*}
\end{prop}

\section[Link with $(-\beta)$-expansions]{Link with expansion in negative bases}

\subsection[$\beta$ attached to several Lyndon words]{Reals associated to several alternate Lyndon words} 

We consider a finite alphabet $ \A $. The aim of this section is to determine for what real number is associated several alternate Lyndon word.

\begin{prop}\label{proposition 9}
Let $\beta \geq 1$ associated to at least two alternate Lyndon words. Then, there exists a weak alternate Lyndon word (a periodic Lyndon word) associated to $ \beta $. 
\end{prop}

\begin{preu}
We consider two alternate Lyndon words over the alphabet $ \A$, $(a_i)_{i \geq 1}$ and $ (d_i)_{i \geq 1}$. Assume that $ (d_i)_{i \geq 1}$ is the alternate Lyndon max associated to $ \beta$. That is,
\begin{equation*}
 a_1a_2a_3 \cdots \prec d_1d_2 d_3 \cdots.
\end{equation*}
This implies that there exists an integer $ n $ such that $ a_i = d_i $ for all 
$ i < n $ and $ (-1)^n(a_n - d_n) < 0$. Furthermore, 

\begin{equation*}
 \sum \limits_{k \geq 1} \frac{a_k}{(-\beta)^k} = \sum \limits_{k \geq 1} \frac{d_k}{(-\beta)^k}.
\end{equation*}
Then :

\begin{equation*}
 a_n - d_n + \sum\limits_{k \geq 1} \frac{a_{n+k}}{(-\beta)^k} - \sum\limits_{k \geq 1} \frac{d_{n+k}}{(-\beta)^k} = 0 .
\end{equation*}
Since $(a_{k+n})_{k \geq 1}$ and $(d_{k+n})_{k \geq 1}$ belong to the Lyndon systems associated to the words $ (a_i)_{i \geq 1} $ and $(d_{i})_{i \geq 1} $; and considering the previous theorem, it follows that :

\begin{equation*}
 \vert \sum\limits_{k \geq 1} \frac{a_{n+k}}{(-\beta)^k} - \sum\limits_{k \geq 1} \frac{d_{n+k}}{(-\beta)^k} \vert \leq 1.
\end{equation*}
But $ \vert a_n - d_n \vert \geq 1 $. Hence :

\begin{minipage}{5.5cm}
\begin{equation*}
 \begin{cases}
  a_n - d_n                &= 1 \\
  \sum\limits_{k \geq 1} \frac{a_{n+k}}{(-\beta)^k} - \sum\limits_{k \geq 1} \frac{d_{n+k}}{(-\beta)^k} &= -1
 \end{cases}
\end{equation*}
\end{minipage} or
\begin{minipage}{5.5cm}
\begin{equation*}
 \begin{cases}
  a_n - d_n &= -1 \\
  \sum\limits_{k \geq 1} \frac{a_{n+k}}{(-\beta)^k} - \sum\limits_{k \geq 1} \frac{d_{n+k}}{(-\beta)^k} &= 1 .
 \end{cases}
\end{equation*}
\end{minipage}

\vspace{1cm}
That is

\begin{minipage}{5.5cm}
\begin{equation} 
\begin{cases} a_n - d_n &= 1  \text{ with $ n $ odd } \\
\sum\limits_{i \geq 1} \frac{a_{n+i}}{(-\beta)^i} &= -\frac{\beta}{\beta+1} \\
 \sum\limits_{i\geq1} \frac{d_{n+i}}{(-\beta)^i} &=  \frac{1}{\beta+1} \label{(24)}
 \end{cases}
\end{equation}
\end{minipage} or
\begin{minipage}{5.5cm}
\begin{equation}
 \begin{cases}
 a_n - d_n &= -1 \text{ with $ n $ even } \\
\sum\limits_{i \geq 1} \frac{a_{n+i}}{(-\beta)^i} &= \frac{1}{\beta+1} \\
\sum\limits_{i \geq 1}\frac{d_{n+i}}{(-\beta)^i}   &= -\frac{\beta}{\beta + 1}.\label{(25)}
 \end{cases}
\end{equation}
\end{minipage}

Let $ u $ be the word such that
\begin{equation*}
 u = \begin{cases}
      \overline{d_1 d_2 \cdots d_n 0} &\text{ if $ n $ is odd }, \\
      \overline{d_1 d_2 \cdots d_n }   &\text{ if $ n $ is even }.
     \end{cases}
\end{equation*}
Since $ (-1)^n(a_n-d_n)<0 $ and for all $ i $, $ 0 \leq d_i \leq d_1$. It follows that $ (a_i)_{i \geq 1} \prec u \preceq (d_i)_{i \geq 1 } $. Moreover, $ u $ is an alternate Lyndon word since $ d_1 \cdots d_n 0 \preceq d_1d_2 \cdots d_nd_{n+1} $ if $ n $ is odd and $ d_1d_2 \cdots d_n d_1 \preceq d_1 \cdots d_{n+1}$ if $ n $ is even. According to the proposition \ref{proposition 2}, $ u $ is attached to $ \beta $. 
Then, $ u $ is a weak alternate Lyndon word associated to $ \beta $ 
\end{preu}

\vspace{1cm}
In the sections \textbf{3.2} and \textbf{3.3}, we saw that the Lyndon words attached to 1 are the terms of the sequence $ (\overline{\phi^n(1)} )_{n \geq 1} $ with the limit included. Moreover, if an alternate Lyndon word $ (d_i)_{i \geq 1} $ satisfied:
\begin{equation*}
 (d_i)_{i \geq 1} \prec \lim \limits_{n \rightarrow + \infty} \phi^n(1)
\end{equation*}
then, it is attached to a real $ \beta $ strictly greater than 1. For such  $ \beta $, we determine the set of the alternate Lyndon which is attached to it.

\begin{prop}\label{proposition 10}
Let $ \overline{ (d_1, d_2, \cdots, d_n) } $ be a weak alternate Lyndon word attached to a real number $ \beta> 1$ (with $ n $ minimal). Then, all alternate Lyndon word such that
 \begin{align}
 d_1d_2 \cdots d_{n-1}(d_n-1) 0 \overline{ d_1 \cdots d_n } \preceq & a_1 a_2 \cdots \preceq \overline{ d_1 d_2 \cdots d_n } \text{ if $ n $ is even } \label{(26)}\\ 
 d_1d_2 \cdots d_n \overline{ d_1 \cdots d_{n-1} (d_n - 1) 0 } \preceq & a_1 a_2 \cdots \preceq \overline{d_1 \cdots d_{n-1}(d_n-1)0 } \text{ if $ n $ odd } \label{(27)}
\end{align}
is attached to $ \beta $.
 \end{prop}

\begin{preu}

We note $ \mathbf{Lyn}(\beta) $ the set of alternate Lyndon words associated to $ \beta $, $ max $ and $ min $ the two words of $ \mathbf{Lyn}(\beta) $ such that for all $ (a_i)_{i \geq 1} \in \mathbf{Lyn}(\beta) $, 
 
 \begin{equation*}
  min \preceq (a_i)_{i \geq 1} \preceq max .
 \end{equation*}
 From the proposition \ref{proposition 9}, we know that there exists an alternate periodic Lyndon word associated to $ \beta$, shall we say $ \overline{( d_1, d_2, \cdots, d_n )} $ (with $ n $ minimal). Then, $ max $ and $ min $ begin by $ d_1 d_2 \cdots d_{n-1} $.
 
 \begin{itemize}
  \item if $ n $ is even 
  \begin{equation*}
   min = d_1 d_2 \cdots d_{n-1} (d_n - 1)0 max 
  \end{equation*}
 since the length of $ d_1 d_2 \cdots d_{n-1} (d_n - 1)0 $ is odd. And
  \begin{equation*}
   max = d_1d_2 \cdots d_n max 
  \end{equation*}
since the length of $ d_1 d_2 \cdots d_n  $ is even. Then, $ max $ is periodic with period $ n $.
 
 \begin{equation*}
  max = \overline{d_1d_2 \cdots d_n }.
 \end{equation*}
And
\begin{equation*}
 min = d_1 d_2 \cdots d_{n-1} (d_n-1) 0 \overline{d_1 d_2 \cdots d_n}.
\end{equation*}

\item If $ n $ is odd, 
\begin{equation*}
 min = d_1 d_2 \cdots d_n max 
\end{equation*}
 since the length of $   d_1 \cdots d_n  $ is odd. And
\begin{equation*}
 max = d_1 d_2 \cdots d_{n-1}(d_n -1)0 max 
\end{equation*}
since the length of $ d_1\cdots d_{n-1}(d_n -1)0 $ is even.
Then $ max $ is periodic with period $ n+1 $.
\begin{equation*}
 max = \overline{ d_1 d_2 \cdots d_{n-1}(d_n - 1) 0}.
\end{equation*}
And
\begin{equation*}
 min = d_1 d_2 \cdots d_n \overline{ d_1 \cdots d_{n-1}(d_n-1) 0 }
\end{equation*}
So, $ \mathbf{Lyn}(\beta) $ is the set of alternate Lyndon words such that

\begin{align*}
 d_1d_2 \cdots d_{n-1}(d_n-1) 0 \overline{ d_1 \cdots d_n } \preceq & a_1 a_2 \cdots \preceq \overline{ d_1 d_2 \cdots d_n } \text{ if $ n $ is even },\\ 
 d_1d_2 \cdots d_n \overline{ d_1 \cdots d_{n-1} (d_n - 1) 0 } \preceq & a_1 a_2 \cdots \preceq \overline{d_1 \cdots d_{n-1}(d_n-1)0 } \text{ if $ n $ odd }.
\end{align*}
\end{itemize}
This completes the proof
\end{preu}

In the previous proposition, we are trying to imply that $ max \prec \lim \limits_{n \rightarrow +\infty} \phi^n (1) $. 

\begin{defi}
In the following, let given a real $ \beta>1 $ associated to an alternate Lyndon word $(d_i)_{i \geq 1}$, we will note $(d^{*}_i)_{i \geq 1 }$ the largest word (in sense of the alternate order) associated to $ \beta $, that is : 

\begin{equation}
(d^*_i)_{i \geq 1} = \begin{cases}
                           \overline{( d_{1},\cdots,d_{2p}, d_{2p+1}-1, 0)} & \text{ if $ \overline{(d_{1}, \cdots, d_{2p+1})}$ is a Lyndon associated to $ \beta $} \\ 
                           (d_{1}, d_{2}, \cdots ) & \text{ otherwise}
                     \end{cases} \label{(28)}
\end{equation}
In the relation above, we consider $ 2p+1$ minimal.
\begin{equation*}
 (d^{*}_i)_{i \geq 1} \prec \lim \limits_{n \rightarrow +\infty} \phi^n(1)
\end{equation*}
\end{defi}
Note that in the first case, we have always $ d_{2p+1} \neq 0$. Indeed, if
 $ \overline{ ( d_1, \cdots, d_{2p}, 0) } $ is an alternate Lyndon word associated to $ \beta $, then the largest (in the sense of the alternate order) Lyndon word attached to $ \beta $ is $ \overline{ (d_1, \cdots, d_{2p-1}, d_{2p} + 1 )} $.  

\begin{prop}\label{proposition 11}
The set of real numbers attached to weak Lyndon word are dense in $ [1, + \infty ) $.
\end{prop}
The previous result is due to the fact that between two alternate Lyndon words, as near as they are, we can find a weak Lyndon word.

\subsection{Expansions in negative base}

The expansion in negative base has been introduced by S. Ito and T. Sadahiro in \cite{MR2534912}. Given $ \beta> 1$, to each real number $ x \in I_{\beta} = [ - \frac{\beta}{\beta + 1 }, \frac{1}{\beta + 1 } ) $ is attached a unique writing obtained by a glouton algorithm. We know that $ -\frac{\beta}{\beta+1} $ is the image by $ f_{\beta} $ of at least one alternate Lyndon word. There exists one of these words which is obtained by such an algorithm. The set of all the expansions in base $-\beta$ of real numbers of the interval $ I_\beta$ is a subset of a Lyndon system associated to $ \beta$. As in the Parry's case (see \cite{MR0142719}), this representation of numbers is extended at all real numbers. 

The $(-\beta)$-transformation is the map $ T_{-\beta}$ from $[l_\beta, r_\beta)$ into itself, defined by

\begin{equation*}
T_{-\beta}(x) = -\beta x - \lfloor - \beta x - l_\beta \rfloor.
\end{equation*}
The $(-\beta)$-expansion of a real $x \in [l_\beta, r_\beta)$ is the writing
\begin{equation*}
d(x, -\beta) = \cdot x_1 x_2 \cdots 
\end{equation*}
with $ x_i = \lfloor -\beta T_{-\beta}^{i-1}(x) - l_\beta \rfloor $. If $ x \not\in [l_\beta, r_\beta) $, one finds the smallest integer $ n$ such that $ \frac{x}{(-\beta)^n} \in I_\beta$. 
\begin{equation*}
d(\frac{x}{(-\beta)^n}, -\beta) = \cdot x_1 x_2 \cdots \Rightarrow d(x, -\beta) = x_1x_2 \cdots x_{n}\cdot x_{n+1} x_{n+2} \cdots .
\end{equation*}
For more details, see \cite{MR2534912}. If there is no ambiguity, we can note $ d(x, -\beta)$ by $ (x_i)_{i \geq 1} $. Set $ d(l_\beta, -\beta) = (d_i)_{i \geq1} $. All $(-\beta)$-expansion $(x_i)_{i \geq 1}$ satisfies:
\begin{equation*}
(d_i)_{i \geq 1} \preceq (x_i)_{i \geq n} \prec (r^*_{i})_{i \geq 1}, \text{ $ \forall n $, }
\end{equation*} 
where
\begin{equation*}
(r^*_i)_{i \geq 1} = \begin{cases}
                           \overline{(0,d_{1},\cdots,d_{2p}, d_{2p+1}-1)} & \text{ if $(d_i)_{i \geq 1}=\overline{(d_{1}, \cdots, d_{2p+1}) }$ } \\ 
                           (0, d_{1}, d_{2}, \cdots ) & \text{ otherwise }.
                     \end{cases} 
\end{equation*}
In particular,
\begin{equation}
(d_i)_{i \geq 1} \preceq (d_i)_{i \geq n} \prec (r^{*}_i)_{i \geq 1} . \label{(29)}
\end{equation}
In fact, $(r^{*}_{i +1})_{i \geq 1}$ is the largest alternate Lyndon word associated to $ \beta $. Then, $(d_i)_{i \geq1}$ is a Lyndon word.

Note $ M$ the Lyndon system  associated to $ (d_i)_{i \geq1}$ (the entropy is $ \log \beta $). In the previous subsection, we saw that the map $f_\beta$ is onto. There exist several words in $M$ which have the same image in the interval $[l_\beta, r_\beta]$. This case occurs for words of the form : 
\begin{equation*}
x_1x_2\cdots x_k y_1y_2\cdots \text{ and } x_1x_2\cdots (x_k+(-1)^k) 0y_1y_2 \cdots 
\end{equation*}
with 
\begin{equation*}
f_{\beta}(y_1y_2y_3 \cdots) = - \frac{\beta}{\beta + 1},
\end{equation*} 
that is
\begin{equation*}
f_{\beta}(0y_1y_2 \cdots ) = \frac{1}{\beta + 1}.
\end{equation*}

The words of $M$ can be seen as representations of real numbers. But, these representations are not unique. The set of $(-\beta)$-expansions corresponds to the subset of $M$ forbidding all word for which the image by $f_\beta$ equals $r_\beta$. 

However, considering a Lyndon word $ (d_i)_{i \geq 1}$, the condition \eqref{(29)} does not imply the existence of the base. In other words, all Lyndon word satisfying \eqref{(29)} is not an expansion of $ l_\beta $ for some $ \beta$.  A counter example is the alternate Lyndon word $ \phi^{\infty}(1)$ or $ \overline{1}$. These two words satisfy \eqref{(29)}. The question is : for what alternate Lyndon word, there exists a base $-\beta$ such that it is the expansion of $ \frac{-\beta}{\beta+1}$?

This question has been broached by Wolfgang Steiner in \cite{MR3028656}. We give another resolution of this subject.

Let $ \mathbf{Lyn} $ be the set of Lyondon word (for the alternate order) $(d_i)_{i \geq 1}$ such that there exists a periodic Lyndon word $\overline{a_1\cdots a_k}$, $ k $ minimal, $ a_k \neq 0 $, $ (d_i)_{i \geq 1} \neq \overline{a_1\cdots a_k}$ and
\begin{equation*}
min_k \preceq (d_i)_{i \geq 1} \preceq max_k \prec \phi^{\infty}(1)
\end{equation*}
with : 
\begin{equation} 
min_k = \begin{cases} 
     a_1 a_2 \cdots a_{k-1}(a_k-1)0 \overline{a_1 a_2 \cdots a_k}    &\text{ if $ k $ even } \\
     a_1 a_2 \cdots a_k \overline{ a_1 a_2 \cdots a_{k-1} (a_k-1)0 } &\text{ if $ k $ odd } \label{(31)}
          \end{cases}
 \end{equation}
\begin{equation}
 max_k = \begin{cases}
    \overline{ a_1 a_2 \cdots a_k }                &\text{ if $ k $ even } \\
    \overline{ a_1 a_2 \cdots a_{k-1} (a_k - 1) 0} &\text{ if $ k $ odd }. \label{(32)}
           \end{cases}
\end{equation}

\begin{equation}
 \mathbf{Lyn} = \{ (d_i)_{i \geq 1} \vert min_k \preceq (d_i)_{i \geq 1} \preceq max_k \prec \phi^{\infty}(1) 
 \text{ and } (d_i)_{i \geq 1} \neq \overline{a_1 \cdots a_k} \} 
 \label{(30)}
\end{equation}
If $ \mathbf{Lyn}$ is empty, we consider that $ k = + \infty $, or the word $(a_i)_{i \geq 1}$ is not periodic.

\begin{theo}\label{theorem 3}
 A sequence $ (d_i)_{i \geq 1} $ on an alphabet $ \A $ is the $ (-\beta)$-expansion of $ -\frac{\beta}{\beta+1} $ for some $ \beta > 1 $ if only if the following conditions are satisfied :
 
 \begin{itemize}
  \item[(a)] for all $ n \geq 1 $, $ (d_i)_{i \geq 1} \preceq (d_i)_{i \geq n} \prec (d^{*}_{i-1})_{i \geq 1} $ , $ d^{*}_0 = 0 $ 
  \item[(b)] $ (d_i)_{i \geq 1} \not\in \mathbf{Lyn} $ 
  \item[(c)] $ (d_i)_{i \geq 1} \prec\phi^{\infty} (1) $
 \end{itemize}

\end{theo}

Before the proof, let clarify a little bit the hypothesis of the theorem. If there exists $ \beta> 1$ such that $ (d_i)_{i \geq 1} $ is the $(-\beta)$-expansion of $-\frac{\beta}{\beta+1}$, then there does not exist an integer $n$ for which 
$ f_{\beta}(d_n d_{n+1} \cdots ) = \frac{1}{\beta + 1} $.

The second point means that if the Lyndon word $ (d_i)_{i \geq 1}$ is associated to a real $ \beta$, which is attached to another alternate Lyndon word, then $(d_i)_{i \geq 1} $ is periodic  (it is the word $ \overline{a_1 a_2 \cdots a_k } $ which does not belong to $ \mathbf{Lyn} $)

The condition (c) traduces the fact that the entropy of the system associated to 
 $ (d_i)_{i \geq 1}$ is positive.

\begin{preu}
\begin{itemize}
\item Let $ \beta> 1$ and $ (d_i)_{i \geq 1} $ the $(-\beta)$-expansion of  $-\frac{\beta}{\beta+1}$ in base $ -\beta$. By definition, the first point of the theorem is automatically satisfied. Furthermore, $ (d_i)_{i \geq 1} $ is an alternate Lyndon word and generates a Lyndon system with an entropy strictly positive. It contains the system associated to $ \lim \limits_{n \rightarrow +\infty} \phi^n(1)=\phi^{\infty}(1) $. This implies the third assertion. Moreover, we saw that $ l_\beta $ has several $(-\beta)$-representations if $(d_i)_{i \geq 1} $ is periodic with odd period. Thereby, if $(d_i)_{i \geq1}$ is not periodic with odd period, $(d_i)_{i \geq 1}$ is the $ max $ of alternate Lyndon words associated to $ \beta$ which is not in $ \mathbf{Lyn}$. The second condition of the theorem is satisfied. On the other hand, if $ (d_i)_{i \geq 1 } $ is periodic with odd period, at the proposition \ref{proposition 10}, we saw that the systems associated to the alternate Lyndon word $(d_i)_{i \geq 1}$ and $(d^{*})_{i \geq 1}$ have the same entropy. Then, 
\begin{equation*}
 (d^{*})_{i \geq 1} \prec \lim \limits_{n \rightarrow + \infty} \phi^{n}(1) 
\end{equation*}
Hence the condition (b).

\item Consider $ (d_i)_{i \geq 1} $ satisfying (a), (b) et (c). From (a), $(d_i)_{i \geq 1}$ is associated to a real $ \beta$ and a dynamical system with entropy $ \log \beta$. From (c), $ \log \beta > 0$ that is, $ \beta> 1$. Let $(a_i)_{i \geq 1}$ the $ (-\beta)$-expansion of $ l_\beta = \frac{-\beta}{\beta+1}$. Thereby,  $ (a_i)_{i \geq 1}$  is an alternate Lyndon word associated to $\beta$. 

If $ \mathbf{Lyn} $ is empty, thus $ (a_i)_{i \geq 1} = (d_i)_{i \geq 1} $ since, in this case, there is one and only one alternate Lyndon word associated to $ \beta$.

If $ \mathbf{Lyn}\neq \varnothing $, $ l_{\beta}$  has several $ (-\beta)$-representations. Then, $ (a_i)_{i \geq 1} $ is periodic. It does not belong to $ \mathbf{Lyn} $. Thus $ (a_i)_{i \geq }=(d_i)_{i \geq 1} $ since, there is one alternate Lyndon word attached to $ \beta $ which is not in $ \mathbf{Lyn} $.
\end{itemize}
This proves the theorem \ref{theorem 3}
\end{preu}

\bibliographystyle{alpha}
\bibliography{Lyndon}
\end{document}